\documentclass[12pt]{article}
\usepackage{amssymb,latexsym}
\usepackage{graphicx}
\usepackage{tikz,pgf}

\title{\textbf{The local $h$-polynomial of the
edgewise subdivision of the simplex}}

\author{Christos~A.~Athanasiadis\\
Department of Mathematics\\
University of Athens\\
Athens 15784, Hellas (Greece)\\
\texttt{caath@math.uoa.gr}}

\date{{\small September 22, 2016}}

  \def\NN{{\mathbb N}}

  \def\ZZ{{\mathbb Z}}
  \def\RR{{\mathbb R}}
  
  \def\raa{{\mathrm {\mathbf a}}}
  \def\ree{{\mathrm {\mathbf e}}}
  \def\ruu{{\mathrm {\mathbf u}}}
  \def\rvv{{\mathrm {\mathbf v}}}

  \def\fF{{\mathcal F}}
  \def\kK{{\mathcal K}}
  
  \def\sS{{\mathcal S}}

  \def\asc{{\rm asc}}
  \def\cone{{\rm cone}}
  
  \def\des{{\rm des}}

  \def\esd{{\rm esd}}
  \def\link{{\rm link}}
  \def\sd{{\rm sd}}

  \def\sm{\smallsetminus}
  \newcommand{\qed}{$\hfill \Box$}

\begin{document}
\maketitle

\newtheorem{theorem}{Theorem}[section]
\newtheorem{proposition}[theorem]{Proposition}
\newtheorem{corollary}[theorem]{Corollary}
\newtheorem{defn}[theorem]{Definition}
\newtheorem{remark}[theorem]{Remark}
\newtheorem{lemma}[theorem]{Lemma}
\newtheorem{example}[theorem]{Example}
\newtheorem{examples}[theorem]{Examples}
\newtheorem{conjecture}[theorem]{Conjecture}
\newtheorem{fact}[theorem]{Fact}
\newtheorem{question}[theorem]{Question}
\newtheorem{observation}[theorem]{Observation}
\newtheorem{claim}[theorem]{Claim}

\begin{abstract}
The $r$-fold edgewise subdivision is a well studied
flag triangulation of the simplex with interesting
algebraic, combinatorial and geometric properties.
An important enumerative invariant, namely the local
$h$-polynomial, of this triangulation is computed
and shown to be $\gamma$-nonnegative by providing
explicit combinatorial interpretations to the
corresponding coefficients. A construction of a flag
triangulation of the seven-dimensional simplex whose
local $h$-polynomial is not real-rooted is also
described.

\end{abstract}

\section{Introduction and results}
\label{sec:intro}

\begin{figure}
\begin{center}
\begin{tikzpicture}[scale=0.9]
\label{fg:esd4}

   \draw(0,0) -- (0,4) -- (4,0) -- (0,0);
   \draw(0,2) -- (2,0) -- (2,2) -- (0,2);
   \draw(0,1) -- (3,1) -- (3,0) -- (0,3) -- (1,3)
              -- (1,0) -- (0,1);

\end{tikzpicture}
\caption{The 4-fold edgewise subdivision of the
2-simplex}
\end{center}
\end{figure}
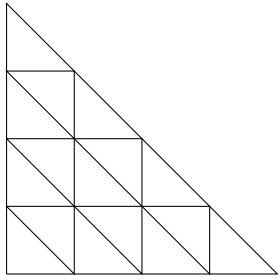

The $r$-fold edgewise subdivision is an elegant
triangulation of a simplicial complex $\Delta$
by which every $k$-dimensional face of $\Delta$
is subdivided into $r^k$ simplices of dimension
$k$. Having arisen in the realm of algebraic
topology~\cite{Fr42}, this construction has
appeared in a wide variety of mathematical
contexts, such as algebraic $K$-theory~\cite{Gr89},
topological cyclic homology~\cite{BHM93}, discrete
and toric geometry~\cite{HPPS, KKMS73},
combinatorial commutative
algebra~\cite{BW09, BR05} and enumerative
combinatorics~\cite{Ath14}. Figure~\ref{fg:esd4}
shows the 4-fold edgewise subdivision of a
two-dimensional simplex. We will denote by
$\esd_r(\Delta)$ the $r$-fold edgewise
subdivision of $\Delta$ and refer to
Section~\ref{sec:sub} for a precise definition.

The edgewise subdivision has been studied in its
own right and shown to have interesting algebraic,
combinatorial and geometric properties; see, for
instance, \cite{BW09, CKW15, EG00, KW12}. Its
effect on the $h$-polynomial, a fundamental
enumerative invariant of a simplicial complex, is
well understood; see Equation~(\ref{eq:hesd}) in
Section~\ref{sec:sub}.

The focus of this paper is on another important
enumerative invariant, namely the local
$h$-polynomial of the $r$-fold edgewise subdivision
of the simplex, which turns out to have a very
elegant combinatorial description. Local
$h$-polynomials were introduced by
Stanley~\cite{Sta92} as a fundamental
tool in his theory of face enumeration for
subdivisions of simplicial complexes. To state the
definition, recall that the $h$-polynomial of a
simplicial complex $\Delta$ of dimension $d-1$ is
defined by the formula

\begin{equation} \label{eq:defh}
h(\Delta, x) \ = \ \sum_{i=0}^d \, f_{i-1} (\Delta)
\, x^i (1-x)^{d-i},
\end{equation}

\noindent
where $f_i (\Delta)$ denotes the number of
$i$-dimensional faces of $\Delta$. Given a
triangulation $\Gamma$ of the abstract simplex $2^V$
on an $n$-element vertex set $V$, the local
$h$-polynomial $\ell_V (\Gamma, x)$ is defined by
the formula

\begin{equation} \label{eq:deflocalh}
  \ell_V (\Gamma, x) \ = \sum_{F \subseteq V} \,
  (-1)^{n - |F|} \, h (\Gamma_F, x),
\end{equation}

\noindent
where $\Gamma_F$ is the restriction of $\Gamma$ to
the face $F \in 2^V$. The importance of local
$h$-polynomials stems from their appearance in the
locality formula \cite[Theorem~3.2]{Sta92}, which
expresses the $h$-polynomial of a triangulation
of a pure simplicial complex $\Delta$ as a sum of
local contributions, one for each face of $\Delta$.
We refer to~\cite{Ath16} for a survey
of basic properties and outstanding open problems
about local $h$-polynomials and for a discussion
of several examples of combinatorial interest.

To state our main results, we need to introduce the
following notation and terminology. We will denote
by $\sS(n, r)$ the set of sequences $w = (w_0,
w_1,\dots,w_n) \in \{0, 1,\dots,r-1\}^{n+1}$ having
no two consecutive entries equal (known as
\emph{Smirnov words}) and satisfying $w_0 = w_n =
0$. We call an index $i \in \{0, 1,\dots,n-1\}$ an \emph{ascent} of such a word $w$ if $w_i < w_{i+1}$
and an index $i \in \{1, 2,\dots,n-1\}$ a
\emph{double ascent} of $w$ if $w_{i-1} < w_i <
w_{i+1}$ (descents and double descents are defined
similarly). We will also denote by ${\rm E}_r$ the
linear operator on the space $\RR[x]$ of polynomials
in $x$ with real coefficients defined by setting
${\rm E}_r (x^m) = x^{m/r}$, if $m$ is divisible
by $r$, and ${\rm E}_r (x^m) = 0$ otherwise.

\begin{theorem} \label{thm:main}
The local $h$-polynomial of the $r$-fold edgewise
subdivision $\esd_r(2^V)$ of the $(n-1)$-dimensional
simplex on the vertex set $V$ can be expressed as

\begin{eqnarray} \label{eq:esdlocal}
\ell_V (\esd_r(2^V), x) &=&
{\rm E}_r \, (x + x^2 + \cdots + x^{r-1})^n \
\ = \ \sum_{w \in \sS(n, r)} x^{\asc (w)},
\end{eqnarray}
where $\asc (w)$ is the number of ascents of $w
\in \sS(n, r)$.
\end{theorem}

An abstract simplicial complex $\Delta$ is called 
\emph{flag} if every clique in the 1-skeleton of 
$\Delta$ (meaning, every set of vertices of $\Delta$ 
pairwise joined by edges) is a face of $\Delta$. 
Barycentric and edgewise subdivisions are
examples of flag triangulations of the simplex. The
following result confirms for edgewise subdivisions
a conjecture of the author~\cite[Conjecture~5.4]{Ath12}
\cite[Conjecture~3.6]{Ath16}, stating that the
polynomial $\ell_V (\Gamma, x)$ is $\gamma$-nonnegative
for every flag triangulation $\Gamma$ of the simplex,
and provides a combinatorial interpretation to the corresponding $\gamma$-coefficients.

\begin{corollary} \label{cor:main}
For all positive integers $n,r$ we have
\begin{eqnarray} \label{eq:edgegamma}
\ell_V (\esd_r(2^V), x) &=&
\sum_{i=0}^{\lfloor n/2 \rfloor} \xi_{n, r, i} \,
x^i (1+x)^{n-2i},
\end{eqnarray}
where $\xi_{n, r, i}$ is the number of sequences
$w = (w_0, w_1,\dots,w_n) \in \sS(n, r)$ with
exactly $i$ ascents which have the following
property: for every double ascent $k$ of $w$ there
exists a double descent $\ell > k$ such that $w_k
= w_\ell$ and $w_k \le w_j$ for all $k < j < \ell$.

In particular, $\ell_V (\esd_r(2^V), -1)$ is equal
to $(-1)^{n/2}$ times the number of sequences $w \in
\sS(n, r)$ with exactly $n/2$ ascents, having this
property.
\end{corollary}

The local $h$-polynomials of flag triangulations
of simplices often have only real roots;
see~\cite[Section~4]{Ath16}. For edgewise
subdivisions, this was shown
recently (using Theorem~\ref{thm:main})
in \cite{Le16a, Le16b} \cite{Zha16}. The following
statement can be viewed as an analogue of a result
of Gal~\cite[Section~3.3]{Ga05}, stating that there
exists a flag triangulation of the five-dimensional
sphere whose $h$-polynomial is not real-rooted.

\begin{theorem} \label{thm:nonreal}
There exists a flag triangulation $\Gamma$ of a
seven-dimensional simplex $2^V$ such that the local
$h$-polynomial $\ell_V(\Gamma, x)$ is not
real-rooted.
\end{theorem}

\medskip
This paper is organized as follows.
Section~\ref{sec:sub} discusses background
on simplicial complexes, subdivisions and local
$h$-polynomials. Theorem~\ref{thm:main} and
Corollary~\ref{cor:main} are proven in
Section~\ref{sec:proofs} using techniques
from enumerative combinatorics. The proof of
Theorem~\ref{thm:nonreal}, given in
Section~\ref{sec:proof}, relies on the
aforementioned construction of Gal and shows
that for every simplicial complex $\Delta$
which is the boundary complex of an
$n$-dimensional convex polytope, there exists
a (regular) triangulation $\Gamma$ (which can
be chosen to be flag, is so is $\Delta$) of a
simplex $2^V$ of dimension $n+1$ such
that $\ell_V (\Gamma, x) = x \, h(\Delta, x)$.
This implies, in particular, that Gal's conjecture
for polytopal flag triangulations of spheres 
follows from the validity 
of~\cite[Conjecture~5.4]{Ath12} for regular flag
triangulations of the simplex. The results of this 
paper were announced without proof 
in~\cite[Section~4]{Ath16}.

\section{Triangulations}
\label{sec:sub}

This section briefly reviews the background on
simplicial and polytopal complexes, their
triangulations and their enumerative invariants
which are needed to understand the main results
and their proofs. For basic notions and more
information on these topics the reader is
referred to the sources~\cite{Bj95, DRS10,
StaCCA, Zi95}. All complexes considered
here are assumed to be finite. The cardinality
and the set of all subsets of a finite set $V$
will denote by $|V|$ and $2^V$, respectively.

\bigskip
\noindent
\textbf{Triangulations.}
Let $\Sigma'$ and $\Sigma$ be geometric simplicial
complexes in some Euclidean space $\RR^N$ (so the
elements of $\Sigma'$ and $\Sigma$ are geometric
simplices in $\RR^N$), with corresponding abstract
simplicial complexes $\Delta'$ and $\Delta$. We say
that $\Sigma'$ is a \emph{triangulation} of $\Sigma$,
and that $\Delta'$ is a \emph{triangulation} of
$\Delta$, if (a) every simplex of $\Sigma'$ is
contained in some simplex of $\Sigma$; and (b) the
union of the simplices of $\Sigma'$ is equal to
the union of the simplices of $\Sigma$. Then, given
any simplex $L \in \Sigma$ with corresponding face
$F \in \Delta$, the subcomplex $\Sigma'_L$ of
$\Sigma'$ consisting of all simplices of $\Sigma'$
contained in $L$ is called the \emph{restriction} of
$\Sigma'$ to $L$. The subcomplex $\Delta'_F$ of
$\Delta'$ corresponding to $\Sigma'_L$ is the \emph{restriction} of $\Delta'$ to $F$. Clearly,
$\Delta'_F$ is a triangulation of the abstract
simplex $2^F$.

Given a simplicial complex $\Delta$ and an
element $v$ not in its vertex set, the \emph{cone}
of $\Delta$ over $v$ is defined as $\cone(\Delta)
= \Delta \cup \{F \cup \{v\}: F \in \Delta\}$.
We note that if $\Delta'$ triangulates $\Delta$,
then $\cone(\Delta')$ naturally triangulates
$\cone(\Delta)$. The \emph{link} of $F \in
\Delta$ is the subcomplex of $\Delta$ defined
as $\link_\Delta (F) = \{ G \sm F: G \in \Delta, \,
F \subseteq G\}$. The \emph{stellar subdivision} of
$\Delta$ on its edge $e = \{a, b\}$ is defined
as the simplicial complex obtained from $\Delta$
by removing all faces which contain $e$ and adding
the sets of the form $F \cup \{v\}$, $F \cup \{v,
a\}$ and $F \cup \{v, b\}$ for all $F \in
\link_\Delta(e)$, where (as before) $v$ is a new
vertex added. This complex is naturally a
triangulation of $\Delta$.

The notions of triangulation of a simplicial
complex and restriction to a face can easily be
extended to polytopal complexes (see
\cite[Section~5.1]{Zi95} for the definition and
examples of polytopal complexes). The boundary
complex of a convex polytope $Q$ (consisting of
all proper faces of $Q$) will be denoted by
$\partial(Q)$.

\bigskip
\noindent
\textbf{Edgewise subdivisions.}
Consider the simplex $2^V$ on the vertex set $V
= \{ \ree_1, \ree_2,\dots,\ree_n \}$ of
coordinate vectors in $\RR^n$. For $\raa = (a_1,
a_2,\dots,a_n) \in \ZZ^n$ set $\iota(\raa) =
(a_1, a_1 + a_2,\dots,a_1 + a_2 + \cdots + a_n)$.
The $r$-fold edgewise subdivision of $2^V$ is the
abstract simplicial complex $\esd_r(2^V)$ on the
vertex set $\Omega_r = \{ (i_1, i_2,\dots,i_n)
\in \NN^n: i_1 + i_2 + \cdots + i_n = r\}$ of
which a set $G \subseteq \Omega_r$ is a face if
$\iota(\ruu) - \iota(\rvv) \in \{0, 1\}^n$, or
$\iota(\rvv) - \iota(\ruu) \in \{0, 1\}^n$, for
all $\ruu, \rvv \in G$.
This is a flag simplicial complex which can be
realized as a triangulation of $2^V$; see, for
instance, \cite[Section~6]{BR05}. To be more
precise, this triangulation is defined by the
dissection of the geometric simplex which is
the convex hull of the set $\Omega_r$ into
smaller geometric simplices by the
affine hyperplanes in $\RR^n$ of the form $x_i
+ x_{i+1} + \cdots + x_j = k$ for $i \le j$
and $k \in \{0, 1,\dots,r\}$. The restriction
of $\esd_r(2^V)$ to a face $F \in 2^V$ coincides
with $\esd_r(2^F)$; it has exactly $r^{\dim(F)}$
faces of the same dimension as $F$. The $r$-fold
edgewise subdivision of an arbitrary simplicial
complex $\Delta$ may be defined so that its
restriction to any face $F \in \Delta$ is
combinatorially isomorphic to $\esd_r(2^F)$;
see~\cite[Section~4]{BW09} \cite[Section~6]{BR05}.

\bigskip
\noindent
\textbf{Face enumeration.}
Let $\Delta$ be a $(d-1)$-dimensional simplicial
complex and let $f_i (\Delta)$ be the number
of $i$-dimensional faces of $\Delta$. The
\emph{$h$-polynomial} of $\Delta$, defined by Equation~(\ref{eq:defh}), is a convenient way to
record the information provided by the numbers $f_i
(\Delta)$. For example, for the two-dimensional
complex of Figure~\ref{fg:esd4} we have $f_{-1}
(\Delta) = 1$, $f_0 (\Delta) = 15$, $f_1 (\Delta)
= 30$, $f_2(\Delta) = 16$ and $h(\Delta, x) =
(1 - x)^3 + 15x(1 - x)^2 + 30x^2(1- x) + 16x^3
= 1 + 12x + 3x^2$. For the importance of
$h$-polynomials, see~\cite[Chapter~II]{StaCCA}.

An explicit formula for the $h$-polynomial of the
$r$-fold edgewise subdivision of $\Delta$ can be
given. Indeed, combining \cite[Corollary~6.8]{BR05}
with \cite[Corollary~1.2]{BW09} (see also
\cite[Equation~(21)]{Ath14}) one gets

\begin{equation} \label{eq:hesd}
h(\esd_r(\Delta), x) \ = \ {\rm E}_r \left(
(1 + x + x^2 + \cdots + x^{r-1})^d \, h (\Delta, x)
\right),
\end{equation}

\medskip
\noindent
where ${\rm E}_r: \RR[x] \to \RR[x]$ is the linear
operator defined in the introduction and $d-1$ is
the dimension of $\Delta$. In particular,

\begin{equation} \label{eq:hesd2V}
h(\esd_r(2^V), x) \ = \ {\rm E}_r \left(
(1 + x + x^2 + \cdots + x^{r-1})^n \right)
\end{equation}

\medskip
\noindent
for every $n$-element set $V$.

Given a triangulation $\Gamma$ of an
$(n-1)$-dimensional simplex $2^V$, the
\emph{local $h$-polynomial}, denoted $\ell_V
(\Gamma, x)$, of $\Gamma$ is defined
\cite[Definition~2.1]{Sta92} by
Equation~(\ref{eq:deflocalh}), where $\Gamma_F$
is the restriction of $\Gamma$ to the face $F \in
2^V$. For the triangulation of Figure~\ref{fg:esd4}
we have $h(\Gamma_V, x) = 1 + 12x + 3x^2$,
$h(\Gamma_F, x) = 1 + 3x$ if $F$ has two elements
and $h(\Gamma_F, x) = 1$ otherwise, so that
$\ell_V(\Gamma, x) = (1 + 12x + 3x^2) - 3(1 + 3x)
+ 3 - 1 = 3x + 3x^2$. The polynomial
$\ell_V(\Gamma, x)$ has degree at most $n-1$
(unless $V = \varnothing$, in which case $\ell_V
(\Gamma, x) = 1$) and nonnegative and symmetric
coefficients, in the sense that $x^n \, \ell_V
(\Gamma, 1/x) = \ell_V (\Gamma, x)$. For examples
and further properties of local $h$-polynomials,
see \cite{Ath16} \cite[Part~I]{Sta92}
\cite[Section~II.10]{StaCCA}.

\bigskip
\noindent
\textbf{Barycentric subdivisions.}
Consider a polytopal complex $\kK$. Choose a
point $p_G$ in the relative interior of each face
$G$ of $\kK$ (and note that the chosen points
include all vertices of $\kK$). The
\emph{barycentric subdivision} of $\kK$ is the
unique triangulation $\sd(\kK)$ of $\kK$ with the
following properties: (a) the vertices of $\sd(\kK)$
are exactly the chosen points, and; (b) the
restriction of $\sd(\kK)$ on any face $G$ of $\kK$
of positive dimension is the cone over $p_G$ of
the restriction of $\sd(\kK)$ on the boundary of
$G$. A more general construction appears in
Section~\ref{sec:proof} (see
Lemma~\ref{lem:partialbary}).

\section{Proof of Theorem~\ref{thm:main} and
Corollary~\ref{cor:main}}
\label{sec:proofs}

This section proves Theorem~\ref{thm:main} and
Corollary~\ref{cor:main} using methods of enumerative
combinatorics. The proof of the latter uses a
variant of the valley hopping technique of
Foata, Sch\"utzenberger and Strehl; see, for instance,
\cite[Section~4.2]{Pet15} and references therein.

\bigskip
\noindent
\emph{Proof of Theorem~\ref{thm:main}.} We set $\Gamma
= \esd_r(2^V)$ and note that the restriction $\Gamma_F$
is the $r$-fold edgewise subdivision $\esd_r(2^F)$ for
every $F \subseteq V$. Thus, by
Equation~(\ref{eq:hesd2V}) we have
\[ h(\Gamma_F, x) \ = \ {\rm E}_r \left(
(1 + x + x^2 + \cdots + x^{r-1})^{|F|} \, \right) \]
for every $F \subseteq V$. The defining
Equation~(\ref{eq:deflocalh}) then yields
\begin{eqnarray*}
\ell_V (\Gamma, x) & = &
{\rm E}_r \left( \, \sum_{k=0}^n \, (-1)^{n-k}
\, {n \choose k}
(1 + x + x^2 + \cdots + x^{r-1})^k \right) \\
& & \\
& = & {\rm E}_r \left( (x + x^2 + \cdots +
x^{r-1})^n \right).
\end{eqnarray*}

\medskip
\noindent
This proves the first equality in (\ref{eq:esdlocal}).
For the second equality, we extend the action of
the operator ${\rm E}_r$ on the space of formal power
series in $x$ with real coefficients in the obvious
way and note that

\begin{eqnarray*}
\frac{h (\Gamma, x)}{(1-x)^n} & = &
\frac{ {\rm E}_r \left(
(1 + x + x^2 + \cdots + x^{r-1})^n \right)}
{(1-x)^n} \ = \
{\rm E}_r \left(
\frac{(1 + x + x^2 + \cdots + x^{r-1})^n}
{(1-x^r)^n} \right) \\
& & \\
& = & {\rm E}_r \left( \frac{1}{(1-x)^n} \right) \
= \ \sum_{m \ge 0} \, {n+rm-1 \choose n-1} \, x^m.
\end{eqnarray*}

We now use the identity
\begin{equation} \label{eq:savage}
\sum_{m \ge 0} \, {n+rm \choose n} \, x^m \ = \
\frac{{\displaystyle \sum_{w \in \{0, 1\dots,r-1\}^n} x^{\asc(w)}}} {(1-x)^{n+1}},
\end{equation}

\medskip
\noindent
where $\asc(w)$ stands for the number of indices $i
\in \{0, 1,\dots,n-1\}$ such that $w_i < w_{i+1}$
for $w = (w_1, w_2,\dots,w_n) \in \{0,
1\dots,r-1\}^n$, with the convention that $w_0 = 0$.
This identity follows from~\cite[Corollary~8]{SS12},
derived in the context of Ehrhart theory, although
a simple combinatorial proof (for instance, one
analogous to the proof of~\cite[Theorem~1]{Pet13})
can be given.

Replacing $n$ with $n-1$ in~(\ref{eq:savage}) and
comparing with our computation of $h(\Gamma, x) /
(1-x)^n$ yields that
\begin{equation} \label{eq:hesdwords}
h(\Gamma, x) \ = \ {\displaystyle
\sum_{w \in \{0, 1\dots,r-1\}^{n-1}} x^{\asc(w)}}.
\end{equation}

\medskip
\noindent
Let us consider the sequences which appear in the
right-hand side of this formula as having length
$n+1$ and first and last coordinate zero and set
$V := \{v_1, v_2\dots,v_n\}$. Then, for every $F
\subseteq V$, we may interpret $h(\Gamma_F, x)$ as
\[ h(\Gamma_F, x) \ = \ {\displaystyle \sum \,
x^{\asc(w)}}, \]

\medskip
\noindent
where the summation runs over all words $w = (w_0,
w_1,\dots,w_n) \in \{0, 1\dots,r-1\}^{n+1}$
satisfying $w_{i-1} = w_i$ for all $i \in \{1,
2,\dots,n\}$ with $v_i \notin F$ and $w_0 = w_n
= 0$. The defining Equation~(\ref{eq:deflocalh})
and an application of the principle of
inclusion-exclusion then yield
\[ \ell_V (\Gamma, x) \ = \sum_{F \subseteq V} \,
(-1)^{n - |F|} \, h (\Gamma_F, x) \ = \
\sum_{w \in \sS(n, r)} x^{\asc (w)} \]

\medskip
\noindent
and the proof follows.
\qed

\bigskip
\noindent
\emph{Proof of Corollary~\ref{cor:main}.} We
define an equivalence relation on the set $\sS(n,r)$
as follows. Let $w = (w_0, w_1,\dots,w_n) \in \sS
(n,r)$ and $k \in \{1, 2,\dots,n-1\}$ be a double
descent of $w$. Since $w_0 = 0$, we may consider the
largest index $1 \le \ell < k$ such that $w_{\ell-1}
< w_k$. We then have $w_\ell \ge w_k$. Assuming
$w_\ell > w_k$, we define the \emph{left match} of $w$
with respect to $k$ as the sequence $w'$ which is 
obtained from $w$ by first deleting $w_k$ and then 
inserting it between $w_{\ell-1}$ and $w_\ell$. 
Formally, we define $w' = (w'_0, w'_1,\dots,w'_n)$ 
by setting $w'_\ell = w_k$, $w'_{i+1} = w_i$ for 
$\ell \le i < k$ and $w'_i = w_i$ for all other values 
of $i$. We note that the left match $w'$ is an 
element of $\sS(n,r)$ with one more ascent than $w$. Similarly,
suppose $k \in \{1, 2,\dots,n-1\}$ is a double ascent
of $w = (w_0, w_1,\dots, w_n) \in \sS(n,r)$. Since
$w_n = 0$, we may consider the smallest index $k <
\ell < n$ such that $w_{\ell+1} < w_k$. We must then
have $w_k \le w_\ell$. Assuming $w_k < w_\ell$, we
define the \emph{right match} of $w$ with respect to 
$k$ as the sequence obtained from $w$ by first 
deleting $w_k$ and then inserting it between $w_\ell$ 
and $w_{\ell+1}$ and note that this 
is an element of $\sS(n,r)$ with one less ascent 
than $w$. For example, if $n=9$, $r = 3$ and $w = 
(0, 2, 1, 2, 1, 0, 1, 2, 1, 0)$, then $k = 4$ is a 
double descent of $w$ and the corresponding left match 
of $w$ is the sequence $(0, 1, 2, 1, 2, 0, 1, 2, 
1, 0)$, obtained from $w$ by deleting its fifth entry 
and inserting it between the first two; see Figure~2. 
Note that if $w'$ is the left match of $w$ with 
respect to its 
double descent $k$, with $w_k$ inserted between 
$w_{\ell-1}$ and $w_\ell$ as above, then the right 
match of $w'$ with respect to its double ascent 
$\ell$ is equal to $w$ and similarly for right 
matches of $w$.

\begin{figure}
\begin{center}
\begin{tikzpicture}[scale=0.6]
\label{fg:class}

\draw(0,0) node(1){};
\draw(2,2) node(2){};
\draw(3,1) node(3){};
\draw(4,2) node(4){};
\draw(5,1) node(5){};
\draw(6,0) node(6){};
\draw(7,1) node(7){};
\draw(8,2) node(8){};
\draw(9,1) node(9){};
\draw(10,0) node(10){};

\draw(1) -- (10);
\draw(1) -- (2) -- (3) -- (4) -- (5) -- (6) -- (7)
-- (8) -- (9) -- (10);
\draw [dashed, ->] (5,1) -- (1,1);

\node[circle,inner sep=0pt,fill=white] at (1) {$\bullet$};
\node[circle,inner sep=0pt,fill=white] at (2) {$\bullet$};
\node[circle,inner sep=0pt,fill=white] at (3) {$\bullet$};
\node[circle,inner sep=0pt,fill=white] at (4) {$\bullet$};
\node[circle,inner sep=0pt,fill=white] at (5) {$\bullet$};
\node[circle,inner sep=0pt,fill=white] at (6) {$\bullet$};
\node[circle,inner sep=0pt,fill=white] at (7) {$\bullet$};
\node[circle,inner sep=0pt,fill=white] at (8) {$\bullet$};
\node[circle,inner sep=0pt,fill=white] at (9) {$\bullet$};
\node[circle,inner sep=0pt,fill=white] at (10) {$\bullet$};

\end{tikzpicture}
\caption{The equivalence class of $(0, 2, 1,
2, 1, 0, 1, 2, 1, 0)$}
\end{center}
\end{figure}
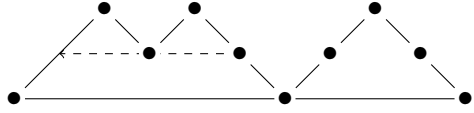

We say that two elements of $\sS(n,r)$ are equivalent
if one can be obtained from the other by a sequence of
matchings, as just described. Our previous discussion 
shows that 
this defines an equivalence relation on $\sS(n,r)$ and 
that the equivalence class of $u \in \sS(n,r)$ has 
exactly $2^{m(u)}$ elements, where $m(u)$ is the total 
number of (left or right) matches of $u$. Moreover, 
each equivalence class contains a unique element with
no right match; these representatives are exactly the 
elements $w \in \sS(n,r)$ which have the property that 
for every double ascent $k$ of $w$ there exists a 
double descent $\ell > k$ such that $w_k = w_\ell$ 
and $w_k \le w_j$ for all $k < j < \ell$. We leave it 
to the reader to verify that $m(w) = \des(w) - \asc(w)
= n - 2\asc(w)$ for every such element $w$. Since 
there are ${m(w) \choose i}$ ways to choose an element 
$u$ in the equivalence 
class $O(w)$ of $w$ by applying $i$ left matchings to 
$w$ and since every such element $u$ has exactly 
$\asc(w) + i$ ascents, we conclude that
\[ \sum_{u \in O(w)} \, x^{\asc(u)} \ = \ 
\sum_{i=0}^{m(w)} {m(w) \choose i} \, x^{\asc(w) + i} \ = 
\ x^{\asc(w)} \, (1+x)^{m(w)} \ = \ x^{\asc(w)}
\, (1+x)^{n - 2\asc(w)}. \]
For our example $w = (0, 2, 1, 2, 1, 0, 1, 2, 1, 0)$,
this expression is equal to $x^4(1+x)$; see
Figure~2 for an attempt to draw the
equivalence class of this sequence. Summing over
all equivalence classes and taking
Theorem~\ref{thm:main} into account, we get the
desired expression for $\ell_V (\esd_r(2^V), x)$.
\qed

\section{Proof of Theorem~\ref{thm:nonreal}}
\label{sec:proof}

This section uses enumerative and geometric
arguments, as well as the construction of flag
triangulations of the five-dimensional sphere whose
$h$-polynomials are not real-rooted by Gal~\cite{Ga05},
to prove Theorem~\ref{thm:nonreal}. Since it will be
crucial in the proof, we note that the spheres
constructed by Gal are easily seen to be polytopal,
meaning they are boundary complexes of simplicial
polytopes. We prepare for the proof with a couple of
lemmas.

\begin{lemma} \label{lem:stellar}
For every triangulation $\Gamma$ of an
$(n-1)$-dimensional simplex $2^V$ having an
interior vertex $p$, there exists a triangulation
$\Gamma'$ of an $n$-dimensional simplex $2^{V'}$
such that $\ell_{V'} (\Gamma', x) = x \,
h(\link_\Gamma(p), x)$. If $\Gamma$ is flag, then
$\Gamma'$ can be chosen to be flag as well.
\end{lemma}

\noindent
\emph{Proof.} Consider the cone of
$\Gamma$ over a new vertex $v$ and let $\Gamma'$ be
the stellar subdivision of this cone on the edge $e
= \{v, p\}$. Then $\Gamma'$ is a triangulation of
the simplex $2^{V'}$, where $V' = V \cup \{v\}$. The
effect of stellar subdivisions on edges on the local
$h$-polynomial was studied in~\cite[Section~6]{Ath12}.
In particular, from the first displayed equation in
the proof of Proposition~6.1 in that reference, the
definition \cite[Equation~(3-5)]{Ath12} of the
relative local $h$-polynomial and the fact that $e$
is an interior face of $\cone(\Gamma)$ (or by direct
computation), we get
\[ \ell_{V'} (\Gamma', x) \ = \
   \ell_{V'} (\cone(\Gamma), x) \, + \,
   x \, h(\link_{\cone(\Gamma)}(e), x) \ = \
   x \, h(\link_\Gamma(p), x). \]
For the second equality we have used that fact
(see the discussion after Proposition~4.14 in
\cite{Sta92}) that the local $h$-polynomial of a
cone vanishes and the obvious equality
$\link_{\cone(\Gamma)}(e) = \link_\Gamma (p)$.
The last statement of the lemma follows from the
previous construction and the fact (see, for
instance, \cite[Proposition~2.4.6]{Ga05}) that
conings and stellar subdivisions on edges preserve
flagness.
\qed

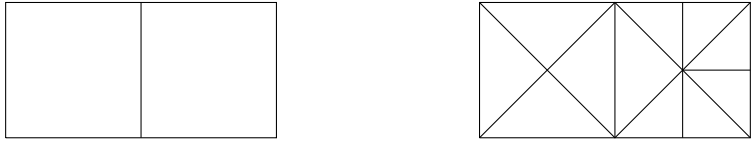
\begin{figure}
\begin{center}
\begin{tikzpicture}[scale=0.9]

\draw(0,0) -- (4,0) -- (4,2) -- (0,2) -- (0,0);
\draw(2,0) -- (2,2);
\draw(7,0) -- (11,0) -- (11,2) -- (7,2) -- (7,0);
\draw(9,0) -- (9,2);
\draw(7,0) -- (9,2) -- (11,0);
\draw(7,2) -- (9,0) -- (11,2);
\draw(10,0) -- (10,2);
\draw(10,1) -- (11,1);

\label{fg:pbsd}
\end{tikzpicture}
\caption{A partial barycentric subdivision}
\end{center}
\end{figure}

\bigskip
The following lemma constructs a partial barycentric subdivision of a polytopal complex $\kK$ with respect
to a simplicial subcomplex $\fF$, which coincides with
the usual barycentric subdivision when $\fF$ consists
only of vertices of $\kK$. Figure~3 shows
this subdivision when $\kK$ consists of two squares
(along with their faces) sharing a common edge
and $\fF$ is the boundary complex of one of them.
This classical construction is usually called the
\emph{barycentric subdivision of $\kK$ relative to
$\fF$}; see, for instance,
\cite[Definition~2.5.7]{Mau80}.

\begin{lemma} \label{lem:partialbary}
Let $\kK$ be a polytopal complex and $\fF$ be a
simplicial subcomplex of $\kK$. Choose a point $p_G$
in the relative interior of each face $G$ of $\kK$,
which is not a face of $\fF$ of positive dimension.
Then there exists a unique triangulation
$\sd_\fF(\kK)$ of $\kK$ satisfying the following
conditions:
\begin{itemize}
\itemsep=0pt
\item[{\rm (i)}]
the vertices of $\sd_\fF(\kK)$ are exactly the
chosen points $p_G$, and

\item[{\rm (ii)}]
the restriction of $\sd_\fF(\kK)$ on any face
$G$ of $\kK$ of positive dimension which is not
a face of $\fF$ is the cone over $p_G$ of the
restriction of $\sd_\fF(\kK)$ on the boundary
of $G$.
\end{itemize}
Moreover, $\fF$ is a subcomplex of $\sd_\fF(\kK)$
and if $\fF$ is flag, then so is $\sd_\fF(\kK)$.
\end{lemma}

\noindent
\emph{Proof.} The existence and uniqueness of
$\sd_\fF (\kK)$ are known and follow by
induction on the dimension of $\kK$.
That $\fF$ is a subcomplex of $\sd_\fF(\kK)$
follows from condition (i). Finally,
assume that $\fF$ is flag. To show that $\sd_\fF
(\kK)$ is flag, we consider a clique $C$ in the
1-skeleton of $\sd_\fF(\kK)$
and verify that $C$ is the vertex set
of some face of $\sd_\fF(\kK)$ as follows. We define
the \emph{rank} of a chosen point $p_G$ as the
dimension of $G$ and proceed by induction on the
maximum rank $m$ of the elements of $C$. Assume
first that $m=0$, meaning that $C$ consists of
vertices of $\kK$. Since two vertices of $\kK$ are
joined by an edge in $\sd_\fF(\kK)$ only if they
are both vertices of $\fF$, the clique $C$
consists of vertices of $\fF$ and therefore forms
the vertex set of a face of $\fF$, hence of
$\sd_\fF(\kK)$ as well, by flagness of $\fF$.
Suppose now that $m \ge 1$ and let $p_G$ be an
element of $C$ of rank $m$. By the construction
of $\sd_\fF(\kK)$, all other elements of $C$ lie
in $G$ and have rank less than $m$. By the induction hypothesis, they form the vertex set of a face of
$\sd_\fF(\kK)$ which is contained in $G$ and hence,
by the construction of $\sd_\fF(\kK)$, the same is
true for the elements of $C$.
\qed

\bigskip
\noindent
\emph{Proof of Theorem~\ref{thm:nonreal}.}
Gal~\cite[Section~3.3]{Ga05} has constructed a
six-dimensional flag simplicial polytope $Q$ for
which the polynomial $h(\partial(Q), x)$ has at
least one non-real root. Thus, in view of
Lemma~\ref{lem:stellar}, it suffices to prove
the existence of a flag triangulation $\Gamma$
of a six-dimensional simplex $\Sigma$ which has an
interior vertex $p$ such that $\link_\Gamma(p)$
is combinatorially isomorphic to $\partial(Q)$.

Consider the hyperplane $H$ in $\RR^7$ consisting
of all points having last coordinate $x_0 = 0$
and the natural projection $\pi: \RR^7 \to H$. Let
$\Sigma$ be a large six-dimensional simplex within
$H$ and let $R$ be a copy of $Q$, embedded in the
hyperplane $x_0 = 1$ of $\RR^7$, which projects
into the interior of $\Sigma$ under the map $\pi$.
Let $P$ be the convex hull of $\Sigma \cup R$.
Then $R$ is a facet of the polytope $P$.
Projecting the faces of $P$ other than $R$ and
$\Sigma$ under $\pi$ we get a polytopal complex $\kK$
which contains $\partial(\pi(R))$ as a subcomplex.
Since the latter is affinely isomorphic to $\partial
(Q)$, and is therefore flag, applying
Lemma~\ref{lem:partialbary} to $\kK$ and $\fF := 
\partial (\pi(R))$ we get a flag triangulation 
$\sd_\fF(\kK)$ of $\kK$ which contains $\partial
(\pi(R))$ as a subcomplex. This triangulation can 
be completed to a flag triangulation $\Gamma$ of 
$\Sigma$ by choosing a point $p$ in the relative 
interior of $\pi(R)$ and defining $\Gamma$ as the 
union of $\sd_\fF(\kK)$ with the set of all cones 
of the faces of $\partial(\pi(R))$ over $p$. Then 
$\link_\Gamma(p) = \partial(\pi(R))$, which is 
affinely isomorphic to $\partial(Q)$, and hence the triangulation $\Gamma$ constructed has the desired 
properties.
\qed

\bigskip
\noindent \textbf{Acknowledgements}. The author
wishes to thank Anders Bj\"orner for suggesting
the Schlegel diagram method, used to construct
the polytopal complex subdividing the simplex
$\Sigma$ in the proof of Theorem~\ref{thm:nonreal},
Francisco Santos for useful comments on the
subdivision of Lemma~\ref{lem:partialbary} and 
an anonymous referee for suggesting improvements
on the presentation.

\end{document}